\def\draft{n}
\documentclass{amsart}
\usepackage[headings]{fullpage}
\usepackage{amssymb,epic,eepic,epsfig}
\usepackage{pb-diagram,lamsarrow,pb-lams}
\usepackage{pstricks}

\theoremstyle{plain}

\newtheorem{theorem}{Theorem}
\newtheorem{proposition}{Proposition}[section]
\newtheorem{lemma}[proposition]{Lemma}
\newtheorem{corollary}[proposition]{Corollary}

\theoremstyle{definition}
\newtheorem{definition}[proposition]{Definition}

\theoremstyle{remark}

\newtheorem{remark}[proposition]{Remark}

\def\printname#1{
        \if\draft y
                \smash{\makebox[0pt]{\hspace{-0.5in}
                        \raisebox{8pt}{\tt\tiny #1}}}
        \fi
}

\newcommand{\psdraw}[2]
         {\begin{array}{c} \hspace{-1.3mm}
        \raisebox{-4pt}{\epsfig{figure=draws/#1.eps,width=#2}}
        \hspace{-1.9mm}\end{array}}

\newlength{\standardunitlength}
\catcode`\@=11
\long\def\@makecaption#1#2{%
     \vskip 10pt

\setbox\@tempboxa\hbox{
       \small\sf{\bfcaptionfont #1. }\ignorespaces #2}%
     \ifdim \wd\@tempboxa >\captionwidth {%
         \rightskip=\@captionmargin\leftskip=\@captionmargin
         \unhbox\@tempboxa\par}%
       \else
         \hbox to\hsize{\hfil\box\@tempboxa\hfil}%
     \fi}
\font\bfcaptionfont=cmssbx10 scaled \magstephalf
\newdimen\@captionmargin\@captionmargin=2\parindent
\newdimen\captionwidth\captionwidth=\hsize
\catcode`\@=12

\def\lbl#1{\label{#1}\printname{#1}}

\def\BZ{\mathbb Z}

\def\BQ{\mathbb Q}
\def\BR{\mathbb R}
\def\BC{\mathbb C}

\def\Ga{\Gamma}

\def\la{\langle}
\def\ra{\rangle}

\def\Ga{\Gamma}

\def\b{\beta}

\def\longto{\longrightarrow}

\def\calB{\mathcal{B}}

\def\Log{\mathrm{Log}}

\def\Li{\mathrm{Li}}
\def\calB{\mathcal{B}}

\def\TB{\mathrm{TB}}

\def\BCs{\mathbb{\BC}^{**}}
\def\hatFT{\widehat{\mathrm{4T}}}
\def\FT{\mathrm{4T}}
\def\wb2C{\widehat{\b_2(\BC)}}
\def\Ker{\mathrm{Ker}}
\def\calD{\mathcal{D}}
\def\hatBC{\widehat{\calB(\BC)}}
\renewcommand\Im{\mathrm{Im}}
\renewcommand\Re{\mathrm{Re}}

\begin{document}


\title[An extended version of additive $K$-theory]{
An extended version of additive $K$-theory}
\author{Stavros Garoufalidis}
\address{School of Mathematics \\
         Georgia Institute of Technology \\
         Atlanta, GA 30332-0160, USA \\ 
         {\tt http://www.math.gatech} \newline {\tt .edu/$\sim$stavros } }
\email{stavros@math.gatech.edu}

\thanks{The author was supported in part by NSF. \\
\newline
1991 {\em Mathematics Classification.} Primary 57N10. Secondary 57M25.
}

\date{April 10, 2007 }


\begin{abstract}
There are two infinitesimal (i.e., additive)
versions of the $K$-theory of a field $F$:
one introduced by Cathelineau, which is an $F$-module, and another
one introduced by Bloch-Esnault, which is an $F^*$-module. Both versions 
are equipped with a regulator map, when $F$ is the field of complex numbers.

In our short paper we will introduce an extended version of Cathelineau's 
group, and a complex-valued regulator map given by the entropy. We will also
give a comparison map between our extended version and Cathelineau's group.

Our results were motivated by two unrelated sources:
Neumann's work on the extended Bloch group
(which is isomorphic to indecomposable $K_3$ of the complex numbers),
and the study of singularities of generating series of hypergeometric
multisums.
\end{abstract}

\maketitle

\tableofcontents

\section{Introduction}
\lbl{sec.intro}

\subsection{Two flavors of  infinitesimal  $K$-theory}
\lbl{sub.additiveK}

There are two infinitesimal (i.e., additive)
versions of the $K$-theory of a field $F$:

\begin{itemize}
\item[(a)]
$\b_n(F)$, introduced by Cathelineau in \cite[Sec.3.1]{Ca2}, 
which is an $F$-module, and 
\item[(b)]
$\TB_n(F)$, introduced by Bloch-Esnault in \cite{BE}, 
which is an $F^*$-module. 
\end{itemize}
Both versions are equipped with a regulator map, when $F=\BC$ is 
the field of complex numbers; see \cite{Ca1} and \cite{BE}.
Let us recall the definition of $\b_2(F)$ and $\TB_2(F)$
for a field $F$ of characteristic zero. Let

\begin{equation}
\lbl{eq.F**}
F^{**}=F\setminus\{0,1\}.
\end{equation}

\begin{definition}
\lbl{def.be}\cite[Thm.3.4]{BE}
For a field $F$ of characteristic zero consider the abelian
group $\TB_2(F)$ which is an $F^*$-module with
generators $\la a \ra$ for $a \in F^{**}$, subject to the 4-term relation:
\begin{equation}
\lbl{eq.4term2}
\la a \ra -\la b \ra + a \star \la \frac{b}{a} \ra +(1-a) \star
\la \frac{1-b}{1-a} \ra=0.
\end{equation}
Here, $\star$ denotes the action of $F^*$ on $\TB_2(F)$.
\end{definition}

The next definition is taken from \cite[Sec.1.1]{Ca1}; see also
\cite[Sec.4.2]{Ca2}.

\begin{definition}
\lbl{def.ca}\cite{Ca1}
If $F$ is a field of characteristic zero, consider the abelian group $\b_2(F)$
which is an $F$-module with generators $\la a \ra$ for $a \in F^{**}$
subject to the 4-term relation:
\begin{equation}
\lbl{eq.4term3}
\la a \ra -\la b \ra + a \la \frac{b}{a} \ra +(1-a)
\la \frac{1-b}{1-a} \ra=0.
\end{equation}
\end{definition}

The groups $\TB_2(F)$ and $\b_2(F)$ fit into a short exact sequence
of abelian groups

\begin{equation}
\lbl{eq.fit}
0 \longto F \longto \TB_2(F) \longto \b_2(F) \longto 0
\end{equation}
where the first and last terms are $F$-modules and the middle is an 
$F^*$-module; see \cite[Eqn.(1.7)]{BE}.

In our short paper we will introduce an extended version $\wb2C$
of Cathelineau's group $\b_2(\BC)$, and a complex-valued regulator map 
given by the entropy. We will also
give a comparison map between our extended version and Cathelineau's group.

Our results were motivated by two unrelated sources:
\begin{itemize}
\item[(a)]
Neumann's work on the {\em extended Bloch group}
and 
\item[(b)]
the study of singularities of generating series of hypergeometric
multisums. 
\end{itemize}

In a sense, the extended Bloch group is forced upon us by the functional
properties (a 5-term relation)
and the analytic continuation of a single special function:
the {\em Rogers dilogarithm}. For a discussion on the extended Bloch group
see \cite{Ne,DZ,GZ}, and for 
its relation to indecomposable $K$-theory $K_3^{\text{ind}}(\BC)$ see
\cite{Ga1}.

In our paper, we will study the functional properties (a 4-term relation)
and the analytic continuation of the entropy function. 
This will naturally lead us
to introduce an extended version of Cathelineau's group.

The 4-term relation for the entropy function has already appeared 
in the context of additive $K$-theory and infinitesimal
polylogarithms; see \cite{Ca1,Ca2,BE} and also \cite[Defn.2.7]{E-VG}.

From our point of view, the entropy function appears in the asymptotics
of classical binomial coefficients, via the Stirling formula. The
relation to singularities of generating series of hypergeometric multisums
will be explained in a separate publication; see \cite{Ga2}.

\subsection{The Stirling formula and the entropy function}
\lbl{sub.stirling}

In this section we will introduce the entropy function and state some
of its elementary properties.

For $a>0$, we can define $(an)!=\Ga(an+1)$, and for $a>b>0$, we can define
as usual $\binom{an}{bn}=(na)!/((nb)!(n(a-b)!)$.

\begin{definition}
\lbl{def.entropy}
Consider the {\em entropy function} $\Phi$, defined by:
\begin{equation}
\lbl{eq.entropy}
\Phi(x)=-x \log (x)-(1-x) \log(1-x).
\end{equation}
for $x \in (0,1)$.
\end{definition}

$\Phi(x)$ is a multivalued function on $\BC\setminus\{0,1\}$, and is given
by the double integral of a rational function as follows from:

\begin{equation}
\lbl{eq.doubleint}
\Phi''(x)=-\frac{1}{x}-\frac{1}{1-x}.
\end{equation}

%

The next lemma links the growth rate of the binomial coefficients with
the entropy function.

\begin{lemma}
\lbl{lem.stirling}
For $a>b>0$, we  have:
\begin{equation}
\lbl{eq.binomiala}
\binom{an}{bn} \sim e^{n a \Phi(\frac{b}{a})} 
\sqrt{\frac{a}{2b(a-b) \pi n}} \left(1+ O\left(\frac{1}{n}\right)\right)
\end{equation}
\end{lemma}

The proof is an application of {\em Stirling's formula}, which computes the 
asymptotic expansion of $n!$ (see \cite{O}):

\begin{equation}
\lbl{eq.stirling}
\log n! \sim n \log n -n + \frac{1}{2} \log n + \frac{1}{2} \log(2 \pi)
+O\left(\frac{1}{n}\right)
\end{equation}

The next lemma gives a 4-term relation for the entropy function.

\begin{lemma}
\lbl{lem.4term}
For $a,b,a+b \in (0,1)$, $\Phi$ satisfies the 4-term relation:
\begin{equation}
\lbl{eq.4term}
\Phi(b)-\Phi(a)+(1-b) \Phi\left(\frac{a}{1-b}\right) -
(1-a) \Phi\left(\frac{b}{1-a}\right)=0.
\end{equation}  
\end{lemma}

The 4-term relation follows from the {\em associativity} of the 
multibinomial coefficients

\begin{equation}
\lbl{eq.assocbinom}
\binom{(a+b+c)n}{an} \binom{(b+c)n}{bn}=
\binom{(a+b+c)n}{bn} \binom{(a+c)n}{an}=
\frac{((a+b+c)n)!}{(an)!(bn)!(cn)!}
\end{equation}
and Lemma \ref{lem.stirling}, and the specialization to $a+b+c=1$. 
In fact, the 4-term relation \eqref{eq.4term} uniquely
determines $\Phi$ up to multiplication by a complex number. See for example,
\cite{Da} and \cite[Sec.5.4,p.66]{AD}.

\subsection{Analytic continuation of the entropy function}
\lbl{sub.anentropy}

In this section we discuss in detail the analytic continuation
of the entropy function, and the 4-term relation. This will be our motivation
for introducing the extended Bloch-Cathelineau group in the next section.
It is clear from the definition that the entropy function is a multivalued
analytic function on the {\em doubly-punctured plane}

\begin{equation}
\lbl{eq.C**}
\BC^{**}=\BC\setminus\{0,1\}.
\end{equation}

Let $\hat{\BC}$ denote the {\em universal abelian cover} of 
$\BC^{**}$. We can represent the Riemann surface $\hat{\BC}$ as follows. 
Let $\BC_{\text{cut}}$ denote $\BC\setminus\{0,1\}$ cut open along each of the 
intervals $(-\infty,0)$ and $(1,\infty)$ so that each real number $r$ 
outside $[0,1]$ occurs
twice in $\BC_{\text{cut}}$. Let us denote the two occurrences of $r$ by 
$r+0i$ and $r-0i$ respectively. It is now easy to see that $\hat{\BC}$
is isomorphic to the surface obtained from $\BC_{\text{cut}} \times 2\BZ \times
2\BZ$ by the following identifications:

\begin{eqnarray*}
(x+0i;2p,2q) & \sim & (x-0i;2p+2,2q) \quad \text{for} \quad x \in (-\infty,0)
\\
(x+0i;2p,2q) & \sim & (x-0i;2p,2q+2) \quad \text{for} \quad x \in (1,\infty)
\end{eqnarray*}
This means that points in $\hat{\BC}$ are of the form $(z,p,q)$ with 
$z \in \BC_{\text{cut}}$ and $p,q$ even integers. 

\begin{definition}
\lbl{def.regulator}
Let us define the entropy function:
\begin{equation}
\lbl{eq.regulator}
\Phi: \hat\BC \longto \BC, \qquad 
\Phi([z;p,q]) =\Phi(z)- \pi i p z +  \pi i q (1-z).  
\end{equation}
\end{definition}

\begin{remark}
\lbl{rem.coordfree}
Since $\hat{\BC}$ is the Riemann surface of the analytic function 
$z\mapsto(\log(z),\log(1-z))$, it follows that a coordinate-free definition
of $\Phi$ can be given by:
\begin{equation}
\lbl{eq.Phi2}
\Phi(\hat{z})= - z \log(\hat{z})-(1-z)\log(1-\hat{z})
\end{equation}
where $\hat{z} \in \hat{\BC}$ and $z$ denotes its projection to $\BC^{**}$.
\end{remark}

Consider the set

\begin{equation}
\lbl{eq.FT}
\FT:=\{(y,x, \frac{x}{1-y}, \frac{y}{1-x}\}
\subset (\BCs)^4
\end{equation}
of 4-tuples involved in the 4-term relation. Also, let

\begin{equation}
\lbl{eq.FT0}
\FT_0:=\{(x_0,\dots,x_3) \in \FT \, | \, 0 < x_0 < x_1 < x_0+x_1 < 1 \}
\end{equation}
and define $\hatFT \subset \hat{\BC}^4 $ to be the component of the
preimage of $\FT$ that contains all points $((x_0;0,0), \dots, (x_4;0,0))$
with $(x_0,\dots,x_4) \in \FT_0$. For a comparison with the 5-term relation,
see \cite[Rem.2.1]{DZ}.

\subsection{An extended version of $\b_2(\BC)$}
\lbl{sub.extended}

We have all the ingredients to introduce 
an extended version of Cathelineau's group $\b_2(\BC)$.

\begin{definition}
\lbl{def.extb}
The {\em extended group} $\widehat{\b_2(\BC)}$ is the $\BC$-module 
generated by the symbols $\la z;p,q \ra$ with $(z;p,q) \in \hat{\BC}$,
subject to the {\em extended 4-term relation}:
\begin{equation}
\lbl{eq.ex4term}
\la x_0;p_0,q_0 \ra 
- \la x_1;p_1,q_1 \ra 
+(1-x_0) \la \frac{x_1}{1-x_0};p_2,q_2 \ra
-(1-x_1) \la  \frac{x_0}{1-x_1};p_3,q_3 \ra =0
\end{equation}
for 
$ ((x_0;p_0,q_0), \dots, (x_3;p_3,q_3)) \in \hatFT$, and
the relations:
\begin{eqnarray}
\lbl{eq.transfer1}
\la x;p,q \ra -\la x;p,q' \ra &=& \la x;p,q-2 \ra- \la x;p,q'-2 \ra 
\\
\lbl{eq.transfer2}
\la x;p,q \ra -\la x;p',q \ra &=& \la x;p-2,q \ra-\la x;p'-2,q \ra 
\end{eqnarray}
for $x \in \BC^{**}$, $p,q,p',q' \in 2 \BZ$.
\end{definition}

It is easy to see that $\Phi$ satisfies the 4-term relations 
\eqref{eq.transfer1} and \eqref{eq.transfer2}.
Moreover, 
since $\hatFT$ is defined by analytic continuation and $\Phi$ satisfies the
4T-relation of Equation \ref{lem.4term}, $\Phi$ satisfies the $\hatFT$
relation. Thus, the following definition makes
sense.

\begin{definition}
\lbl{def.regulator2}
$\Phi$ gives rise to a 
{\em regulator map}:
\begin{equation}
\lbl{eq.R}
R: \widehat{\b_2(\BC)} \longto \BC
\end{equation}
\end{definition}

Our main theorem compares our extended group $\wb2C$ with Cathelineau's
group $\b_2(\BC)$.

\begin{theorem}
\lbl{thm.1}
\rm{(a)} There is a well-defined map:
\begin{equation}
\lbl{eq.pi}
\pi: \wb2C \longto \b_2(\BC), \qquad \pi(\la x;p,q \ra)=\la x \ra.
\end{equation}
\newline
\rm{(b)}
We have a short exact sequence of $\BC$-modules:
\begin{equation}
\lbl{eq.compare}
0 \longto \BC \stackrel{\chi}\longto \wb2C \stackrel{\pi}
\longto \b_2(\BC) \longto 0
\end{equation}
where 
\begin{equation}
\lbl{eq.chi}
\chi(z)=z \cdot c
\end{equation}
and 
$c= \la x;2,-2 \ra - \la x;0,0 \ra$ is independent of $x \in \BC^{**}$.
\newline
\rm{(c)} The composition $R \circ \chi$ is given by:
\begin{equation}
\lbl{eq.Rz}
(R \circ \chi)(z)=-2 \pi i z.
\end{equation}
In particular, the restriction of $R$ on $\Ker(\pi)$ is 1-1.
\end{theorem}

An important ingredient for the proof of Theorem \ref{thm.1} 
is the following description of the $\hatFT$ relation. 

\begin{definition}
\lbl{def.VFT}
Let $V \subset (\BZ \times \BZ)^4$ be the subspace
\begin{equation}
\lbl{eq.V}
V:=\{((p_0,q_0),(p_1,q_1),(-q_0+p_1,-q_0+r),(p_0-q_1,r-q_1)) \,| \,
p_0,p_1,q_0,q_1,r \in 2 \BZ\}.
\end{equation}
\end{definition}

Let us define 

\begin{equation}
\lbl{eq.4T+}
\FT^+:=\{(x_0,\dots,x_3) \in \FT \, | \,  \Im(x_i)>0, \, i=0,\dots,3 \}
\end{equation}
and let $\hatFT^+$ denote the unique component
of the inverse image of $\FT$ in $\hat{\BC}^4$ which includes the points
$((x_0;0,0),\dots, (x_3;0,0))$ with $(x_0,\dots,x_3) \in \hatFT^+$.

\begin{proposition}
\lbl{prop.1}
We have:
\begin{equation}
\lbl{eq.FTFT}
\hatFT=\hatFT^+ + V=\{\mathbf{x}+\mathbf{v} \, | \, \mathbf{x} \in \hatFT^+
\, \mathrm{and} \, \mathbf{v} \in V \},
\end{equation}
where we are using addition to denote the action of $(\BZ \times \BZ)^4$
by covering transformations on $\hat{\BC}^4$.
\end{proposition}

\begin{corollary}
\lbl{cor.1}
For $(y,x,x/(1-y),y/(1-x)) \in \hatFT^+$ and $p_0,q_0,p_1,q_1,r \in 2 \BZ$
we have:
\begin{equation}
\lbl{eq.exx}
\la y;p_0,q_0 \ra 
- \la x;p_1,q_1 \ra 
+(1-y) \la \frac{x}{1-y}; q_0+p_1, -q_0+r \ra
-(1-x) \la  \frac{y}{1-x};p_0+q_1,r-q_1 \ra =0.
\end{equation}
\end{corollary}

\subsection{A comparison of regulators}
\lbl{sub.regulators}

In \cite[Prop.7]{Ca2} Cathelineau defined purely real/imaginary regulators
on the $\BC$-modules $\b_n(\BC)$ for all $n$:
\begin{equation}
\lbl{eq.cat1}
d\tilde{\calD}_n: \b_n(\BC) \longto \BR(n-1)
\end{equation}
where, we follow the notation of \cite{Ca2} and, as usual, 
$A(n)=(2 \pi i)^n A$. The regulators $d\tilde{\calD}_n$ also known by the name
of {\em infinitesimal polylogarithms} can be
thought as an infinitesimal analogue of the 
{\em modified polylogarithms} $\calD_n$:
\begin{equation}
\lbl{eq.calDn}
\calD_n: \BC \longto \BR(n-1), \qquad
\calD_n(z)=\Re/\Im_n \left(\sum_{k=0}^{n-1} \frac{B_{k} 2^k}{k!} 
\log^k|z| \Li_{n-k}(z)
\right)
\end{equation}
where $\Re/\Im_n$ denotes the real (resp. $i$ times the imaginary) part if 
$n$ is odd (resp. even), $B_k$ are the Bernoulli numbers, and $\Li_k(z)=\sum_k 
z^k/k^n$ is the classical polylogarithm. 
The modified polylogarithms were introduced
by Zagier in \cite{Za}
(generalizing the $n=2$ case  of the $n=2$ Bloch-Wigner 
polylogarithm) and have a motivic
interpretation as a reality condition of Mixed Hodge Structures; see 
\cite{BD}.

The infinitesimal polylogarithms $d\tilde{\calD}_n$ 
have also been studied by Elbaz-Vincent and Gangl; see \cite[Defn.2.7]{E-VG}.
For $n=2$, we have, up to a sign:

\begin{equation}
\lbl{eq.cat2}
i r: \b_2(\BC) \longto i\BR
\end{equation}
where

\begin{equation}
\lbl{eq.r}
r(x)=-x \log|x|-(1-x)\log|1-x|
\end{equation}
is a modification of the entropy function. On the other hand, we have
a complex-valued regulator map $R$ on $\wb2C$. We caution that 
$i \Re \circ R \neq i r \circ \pi$.

Nevertheless, let $\widehat{\b_2(\BR)}$ denote the $\BR$-subspace of 
$\b_2(\BC)$
generated by the real points in $\hat{\BC}$; i.e., the points of the
form $(x;p,q)$ with $x \in \BR^{**}$. The following corollary follows
immediately from the definitions.

\begin{corollary}
\lbl{cor.reg}
The following diagram commutes:
$$
\divide\dgARROWLENGTH by2
\begin{diagram}
\node{\widehat{\b_2(\BR)}}
\arrow{e,t}{R}
\arrow{s,l}{\pi}
\node{\BC}
\arrow{s,r}{i \Re}
\\
\node{\b_2(\BR)}
\arrow{e,t}{i r}
\node{i \BR}
\end{diagram}
$$
\end{corollary}

\subsection{A comparison table}
\lbl{sub.table}

As we mentioned before, there is a close analogy between our paper and
Neumann's work on the extended Bloch group $\hatBC$; see \cite{Ne,DZ,GZ,Ga1}.
Let us summarize this into a table, which enhances the table of 
\cite[Sec.5]{BE}.

\begin{center}
\begin{tabular}{|l|l|l|}
\hline
Theory & $K$-theory & Infinitesimal $K$-theory \\ \hline
Group & $\calB(\BC)$ & $\b_2(\BC)$  \\ \hline
Extended group & $\hatBC$ & $\wb2C$ \\ \hline
Regulator & Rogers dilogarithm & Entropy \\ \hline
 & $R: \hat{\BC} \to \BC/\BZ(2)$ & $R: \hat{\BC} \to \BC$ \\ \hline
Functional Equation & 5-term relation & 4-term relation \\ \hline
Exact sequence & $ 0 \to \BQ/\BZ \to \hatBC \to \calB(\BC) \to 0$
& $ 0 \to \BC \to \wb2C \to \b_2(\BC) \to 0 $ \\ \hline
\end{tabular}
\end{center}

\subsection{Acknowledgement}
An early version of this paper was presented in a talk in the University
of Chicago in the spring of 2007. The author wishes to thank M. Kontsevich, 
for stimulating conversations, C. Zickert for pointing out a sign error
in an earlier version of the paper, and S. Bloch for his
hospitality and encouragement.

\section{Proof of Proposition \ref{prop.1}}
\lbl{sec.prop1}

In this section we will give a proof of Proposition \ref{prop.1}, following
the idea from the Appendix of \cite{DZ}.

Recall the definition of $\FT^+$ from Equation \eqref{eq.4T+}. 
Fix a point $(x_0,\dots,x_3)$ and let

\begin{equation}
P=((x_0;0,0), \dots, (x_3;0,0)) \in \FT^+. 
\end{equation}
Let $x_0=y$ and $x_1=x$. Consider the curve in $\hatFT$ obtained by keeping
$x_1$ fixed and letting $x_0$ move along a closed curve in 
$\BC\setminus\{0,1,1-x_1\}$. By an analysis of the 4-term relation we can 
examine how the values of the $p_i$'s, $q_i$'s change when $x_0$ moves around.
This is recorded in Figure \ref{f.motion1}. It follows that if $x_0$ traverses
a closed curve going $p_0/2$ times clockwise around $0$, 
followed by $q_0/2$
times clockwise around $1$, followed by $r/2$ times clockwise
around $1-x_1$, for $p_0,q_0,r$ even, then the curve in $\hatFT$ ends in $P'$, 
where 

\begin{equation*}
P'= ((x_0;p_0,q_0), (x_1;0,0), (x_2,-q_0,q_0+r), (x_3;p_0,r)).
\end{equation*}

Starting at $P'$, follow the curve in $\hatFT$ obtained by keeping
$x_0$ fixed and letting $x_1$ move a curve going $p_1/2$ times clockwise
around $0$ followed by $q_1/2$ times clockwise around $1$ for
$p_1, q_1$ even. Then the curve in $\hatFT$ ends in $P''$, 
where 

\begin{equation*}
((x_0;p_0,q_0), (x_1;p_1,q_1), (x_2; -q_0+p_1, -q_0+r),
(x_3;p_0+q_1,r-q_1).
\end{equation*}
Since we can connect a point $P$ in $\hatFT$ to a point in the preimage of
$\FT_0$ by first sliding $x_0$ into
the interval $(0,1)$ and then sliding $x_1$ into the interval $(x_0,1)$,
Proposition \ref{prop.1} follows.

\begin{figure}[htpb]
$$
\psdraw{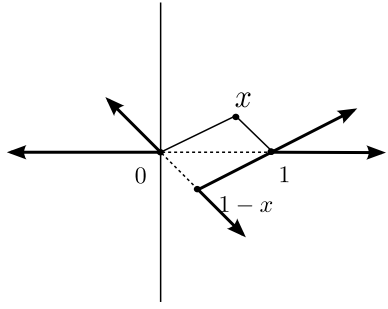}{2in}
$$
\caption{The 6 solid half-lines in the figure are the cuts of the function
$z \mapsto (\Log(z),\Log(1-z))$ in the $x_i$-plane for $i=0,2,3$, when 
$x=x_1$ is fixed. The relevant values of $p_i$ and $q_i$ increase by $2$ 
whenever $y=x_0$ crosses the revelant line in the direction indicated by
the arrows.}\lbl{f.motion1}
\end{figure}

\section{Proof of Theorem \ref{thm.1}}
\lbl{sec.thm1}

For part (a) of Theorem \ref{thm.1}, we need to show that the image 
under $\pi$ of the
extended 4-term relation of Corollary \ref{cor.1} is zero 
in  $\b_2(\BC)$.
Let us replace $x$ by $1-x$. It suffices to show that for all $x,y \in 
\BC^{**}$
with $x \neq y$ we have:

\begin{equation}
\lbl{eq.c1}
\la y \ra -\la 1-x \ra + (1-y)  \la \frac{1-x}{1-y} \ra 
-x\la \frac{y}{x} \ra =0.
\end{equation} 

In \cite[Sec.1.1]{Ca1} Cathelineau proves that the 4-term relation 
\eqref{eq.4term2} implies the following relations:

\begin{eqnarray*}
\la a \ra &=& \la 1-a \ra \\
\la 1/a \ra &=& -1/a \la a \ra
\end{eqnarray*}

It follows that

\begin{eqnarray*}
\la y \ra -\la 1-x \ra + (1-y)  \la \frac{1-x}{1-y} \ra 
-x\la \frac{y}{x} \ra &=& 
\la y \ra -\la 1-x \ra - (1-x)  \la \frac{1-y}{1-x} \ra 
-x\la \frac{y}{x} \ra \\
&=&
\la y \ra -\la x \ra - (1-x)  \la \frac{1-y}{1-x} \ra 
-x\la \frac{y}{x} \ra \\
&=& 0
\end{eqnarray*}
by \eqref{eq.4term3} for $x=a$, $y=b$.

In the remaining of the section we will prove the other parts
of Theorem \ref{thm.1}, following arguments similar to \cite[Sec.7]{Ne}.

\begin{lemma}
\lbl{lem.1}
For $x,y \in \BC^{**}$, $p_0,q_0,p_1,q_1 \in 2 \BZ$ we have:
\begin{equation}
\lbl{eq.2t1}
\la y;p_0-2,q_0+2 \ra - \la y;p_0,q_0 \ra=
\la x;p_1-2,q_1+2 \ra -\la x;p_1,q_1 \ra
\end{equation}
\end{lemma}

\begin{proof}
Replace $(p_0,p_1,q_0,q_1,r)$ by $(p_0-2,p_1-2,q_0+2,q_1+2,r+2)$ in the
4-term relation of Proposition \ref{prop.1} and subtract. Of the 8 terms,
4 cancel and the remaining 4 give the identity \eqref{eq.2t1}.
\end{proof}

For every $x \in \BC^{**}$, let us define

\begin{equation}
\lbl{eq.cx}
c_x=\la x;2,-2 \ra - \la x;0,0 \ra \in \wb2C.
\end{equation}

Setting $(p_0,q_0,p_1,q_1)=(2,-2,2,-2)$
in \eqref{eq.2t1} and recalling \eqref{eq.cx}
implies that 

\begin{equation}
\lbl{eq.c}
c:=c_x
\end{equation}
is independent of $x \in \BC^{**}$.
In addition, the left hand side of \eqref{eq.2t1} is independent of $p_1$
and $q_1$ thus so is the right hand side. Thus, 
\begin{equation}
\lbl{eq.indep}
\la x;p_1-2,q_1+2 \ra-\la x;p_1,q_1 \ra
\end{equation}
is independent of $p_1$ and $q_1$.

\begin{lemma}
\lbl{lem.2}
For all $x \in \BC^{**}$, $p,q \in 2 \BZ$ we have:
\begin{equation}
\lbl{eq.2t2}
\la x; p,q \ra=\frac{1}{4} 
\left( pq \la x; 2,2 \ra - p(q-2) \la x; 2,0 \ra
- q(p-2)\la x; 0,2 \ra + (pq-2p-2q+4)\la x; 0,0 \ra \right).
\end{equation}
\end{lemma}

\begin{proof}
Equation \eqref{eq.transfer1} implies that for all $p \in 2\BZ$ we have:
\begin{equation}
\lbl{eq.nep}
\la x;p,q \ra= \frac{q}{2} \la x;p,2 \ra - \frac{q-2}{2} \la x;p,0 \ra.
\end{equation}
Equation \eqref{eq.transfer2} implies that for all $q \in 2 \BZ$ we have:
\begin{equation}
\lbl{eq.neq}
\la x;p,q \ra= \frac{p}{2} \la x;2,q \ra - \frac{p-2}{2} \la x;0,q \ra.
\end{equation}
Using \eqref{eq.neq}, and telescoping, we can expand the right hand side 
of \eqref{eq.nep}. The result follows.
\end{proof}

Using \eqref{eq.2t2}, it follows that both sides of \eqref{eq.2t1}
are linear functions on $p_0,q_0,p_1,q_1$. In particular, the coefficient
of $p_1$ vanishes. This implies that for all $x \in \BC^{**}$ we have:

\begin{equation}
\lbl{eq.2t3}
\la x;2,2 \ra-\la x;2,0 \ra-\la x;0,2 \ra+\la x;0,0 \ra=0.
\end{equation}

Let $K=\Ker(\pi)$ denote the kernel of $\pi:\wb2C \to \b_2(\BC)$.
It follows that $K$ is the $\BC$-span of 
\begin{equation*}
\la x;p+2,q \ra-\la x;p,q \ra , \qquad \la x;p,q+2 \ra-\la x;p,q \ra
\end{equation*}
for $x \in \BC^{**}$, $p,q \in 2 \BZ$.
Using Lemma \ref{lem.2}, it follows that $K$ is the $\BC$-span of
\begin{equation*}
\la x;2,0 \ra-\la x;0,0 \ra , \qquad \la x;0,2 \ra-\la x;0,0 \ra,
\qquad \la x;2,2 \ra-\la x;0,0 \ra.
\end{equation*}
Using Equation \eqref{eq.2t3}, it follows that 
$K$ is the $\BC$-span of $\la x;0,2 \ra-\la x;0,0 \ra$ 
(for all $x \in \BC^{**}$) and $c$.

For $x \in \BC^{**}$, let us denote

\begin{equation}
\lbl{eq.xbr}
\{x\}:=-\frac{1}{1-x} \left(\la x;0,2 \ra-\la x;0,0 \ra \right).
\end{equation}

Set $p_0=q_0=p_1=q_1=0$
in the 4-term relation \eqref{eq.exx}. It follows that

\begin{equation}
\lbl{eq.2t4}
\la y;0,0 \ra - \la x;0,0 \ra +(1-y)
\la \frac{x}{1-y};0,r \ra -(1-x) \la \frac{y}{1-x};0,r \ra =0.
\end{equation}

Replace $r$ by $r+2$, subtract and divide by $1-x-y$. It follows that

\begin{equation}
\lbl{eq.2t6}
\{ \frac{x}{1-y} \}= \{ \frac{y}{1-x} \}
\end{equation}
which implies that for all $x,y \in \BC^{**}$, we have:
\begin{equation}
\lbl{eq.2t7}
\{ x \}= \{ y \}.
\end{equation}
Thus, $K$ is the $\BC$-span of $\{x\}$ and
$c$, where both elements are independent of $x$.

Going back to the 4-term relation \eqref{eq.exx}, 
replace $q_1$ with $q_1+2$ and subtract. We get the 4-term relation:

\begin{eqnarray*}
\la x; p_1, q_1 \ra - \la x; p_1, q_1+2 \ra &=&
(1-x) \left( \la \frac{y}{1-x}; p_0+q_1+2, r-q_1-2 \ra -
\la \frac{y}{1-x}; p_0+q_1, r-q_1 \ra \right) \\
&=& (1-x) \left( \la \frac{y}{1-x}; 2, -2 \ra -
\la \frac{y}{1-x}; 0, 0 \ra \right) \\
&=&
(1-x) c.
\end{eqnarray*}
Now set $p_1=0$ and $q_1=2$. 
It follows that $\{x\}=c$. Thus, $K$ is the $\BC$-span
of $\{x\}=c$ which is independent of $x$. This implies that $\chi$ is 
well-defined and that $\text{Image}(\chi)=\Ker(\pi)$. 
Since $(R\circ \chi)(z)=z R(c_x)=-2 \pi i z$, it follows that $R \circ
\chi$, and therefore $\chi$, is 1-1.
This concludes the proof of Theorem \ref{thm.1}.

\ifx\undefined\bysame
        \newcommand{\bysame}{\leavevmode\hbox
to3em{\hrulefill}\,}
\fi


\begin{thebibliography}{[EMSS]}


\bibitem[AD]{AD} J. Acz\'el and J. Dhombres, 
        {\em Functional equations in several variables},
        Encyclopedia of Mathematics and its Applications {\bf 31} 
        Cambridge University Press, Cambridge, 1989.

\bibitem[BD]{BD} A. Beilinson and P. Deligne, 
        {\em Interpr\'etation motivique de la conjecture de Zagier reliant 
        polylogarithmes et r\'egulateurs},
        in Motives, Proc. Sympos. Pure Math., {\bf 55}, Part 2, (1994) 97--121.

\bibitem[BE]{BE} S. Bloch and H. Esnault, 
        {\em The additive dilogarithm}, in Kazuya Kato's fiftieth birthday.  
        Doc. Math.  (2003)  Extra Vol. 131--155.
 
\bibitem[Ca1]{Ca1} J.L. Cathelineau,
        {\em  Remarques sur les différentielles des polylogarithmes 
        uniformes},
        Ann. Inst. Fourier {\bf 46}  (1996) 1327--1347. 

\bibitem[Ca2]{Ca2} \bysame,
        {\em Infinitesimal polylogarithms, multiplicative presentations
        of Kaehler differentials and Goncharov complexes}, talk at the\
        workshop on Polylogarithms, Essen, May 1-4.

\bibitem[Da]{Da} Z. Dar\'oczy, 
        {\em Generalized information functions},  
        Information and Control  {\bf 16}  (1970) 36--51.

\bibitem[DZ]{DZ} J.L. Dupont and C. Zickert, 
        {\em A dilogarithmic formula for the Cheeger-Chern-Simons class},
        Geom. Topol. {\bf 10} (2006) 1347--1372. 
 
\bibitem[E-VG]{E-VG} P. Elbaz-Vincent and H. Gangl,
        {\em On poly(ana)logs. I},
        Compositio Math.  {\bf 130}  (2002) 161--210.

\bibitem[Ga1]{Ga1} S. Garoufalidis,
        {\em $q$-terms, singularities and the extended Bloch group},
        preprint 2007 {\tt arXiv:0708.0018}.

\bibitem[Ga2]{Ga2} \bysame,
        {\em An ansatz for the singularities of hypergeometric multisums},
        preprint 2007 {\tt arXiv:0706.0722}.

\bibitem[GZ]{GZ} S. Goette and C. Zickert,
        {\em The Extended Bloch Group and the Cheeger-Chern-Simons Class},
        Geom. Topol. {\bf 11} (2007) 1623--1635.

\bibitem[O]{O} F. Olver, 
        {\em Asymptotics and special functions}, Reprint. 
        AKP Classics. A K Peters, Ltd., Wellesley, MA, 1997.

\bibitem[Ne]{Ne} W.D. Neumann,
        {\em Extended Bloch group and the Cheeger-Chern-Simons class},  
        Geom. Topol. {\bf 8} (2004) 413--474.

\bibitem[Za]{Za} D. Zagier,
        {\em Polylogarithms, Dedekind zeta functions and the algebraic 
        $K$-theory of fields}, 
        in Arithmetic algebraic geometry,
        Progr. Math., {\bf 89} (1991) 391--430.

\end{thebibliography}
\end{document}